# An observation concerning uniquely ergodic vector fields on 3-manifolds

Clifford Henry Taubes[†]

Department of Mathematics
Harvard University
Cambridge, MA 02138

[†] Supported in part by the National Science Foundation

# 1. Introduction

This note concerns the dynamics generated by a certain sort of divergence free vector field on a compact 3-manifold with a prescribed volume form. Let M denote the 3-manifold in question, let $\Omega$ denote the volume 3-form and let v denote the vector field. Use $\omega_v$ to denote the 2-form $\Omega(v, \cdot)$. This is a closed 2-form. Say that v is *exact* when $\omega_v$ has an anti-derivative. This is to say that there exists a 1-form $\upsilon$ such that $d\upsilon = \omega_v$. When v is exact, choose such an anti-derivative so as to define

$$s_v = \int_M \upsilon \wedge \omega \;.$$

(1.1)

This integral does not depend on the chosen anti-derivative.

Motion along v's integral curves defines a 1-parameter family of diffeomorphism $\phi: \mathbb{R} \times M \to M$. A measure, $\sigma$, on M is said to be *v-invariant* when $\phi(t, \cdot)_*\sigma = \sigma$ for all t $\in$ M. By assumption, integration with respect to the volume form on M is a v-invariant measure. The vector field v is said to be *uniquely ergodic* if the volume form is the only v-invariant measure.

What follows is the main result of this note:

**Theorem 1**: *Let* v *denote an exact, uniquely ergodic vector field on* M. *Then* $s_v = 0$.

Truth be told, this theorem may be of little interest because only few vector fields are known to be uniquely ergodic. For example, a vector field with a zero, or one with a closed integral curve is not uniquely ergodic. In any event, an example is described in the final section of this note.

This theorem is proved using the Seiberg-Witten equations in a manner much like that used by the author in [T1], [T2] to prove that the Reeb vector field for a contact 1-form on a 3-manifold has at least one closed orbit. This is the case where $\omega_v$ can be written as da where a has positive pairing with v.

The Seiberg-Witten equations can say more about the dynamics if more is assumed about v. The case when v is the Reeb vector field of a contact 1-form is discussed in [T1], [T2]; and much more is said in [T3] and its sequels. A stable Hamiltonian structure generalizes the notion of a contact structure. The vector field v is the Reeb vector field for a stable Hamiltonian structure in the case where $\omega_v = f$ da where $f$ is a function on the manifold and the 1-form a has positive pairing with v. The Seiberg-Witten equations are used in [HT] to prove that the Reeb vector field for a stable Hamiltonian structure has a closed orbit unless the 3-manifold is a torus bundle over a circle. (This same result was proved more or less at the same time using different techniques by [R].)

A slight generalization of the contact structure case posits $\omega_v = da$ and $a(v) \geq 0$. The Seiberg-Witten equations say the following about this last case:

**Theorem 2**: *Suppose that* $\Omega(v, \cdot) = \omega_v$ *can be written as* $da$ *such that* $a(v) \geq 0$. *Either* v *has a closed integral curve or there is a* v*-invariant measure with support in the set where* $a(v) = 0$.

This theorem says nothing if $a(v)$ is everywhere zero.
Theorems 1 and 2 are proved in the upcoming Section 5.

## 2. The Seiberg-Witten equations

If v has a zero, then v can not be uniquely ergodic. This understood, assume in what follows that v is nowhere zero. Fix a smooth, Riemannian metric on M whose volume form is $\Omega$ and is such that $|v| = 1$. Orient M using $\Omega$ so as to define the principle SO(3) bundle of orthonormal frames in TM. Choose a $Spin_{\mathbb{C}}$ structure $F \to M$, this a U(2) lift of this principle SO(3) frame bundle. With F chosen, introduce the associated vector bundles $S_F = F \times_{U(2)} U(1)$ and let $\mathbb{S} = F \times_{U(2)} \mathbb{C}^2$ where U(2) acts on U(1) via the determinant representation, and where it acts on $\mathbb{C}^2$ via the defining representation.

Fix a smooth 1-form $\mu$ with $C^3$ norm bounded by 1 and a number $r \geq 1$ so as to defined a certain version of the Seiberg-Witten equations. The version in question asks that pair $(\mathbb{A}, \psi)$ where of unitary connection $\mathbb{A}$ on $S_F$ and section $\psi$ of $\mathbb{S}$ obey

$$*F_{\mathbb{A}} = 2r(\psi^\dagger \tau \psi - i*\omega) - i*d\mu \quad and \quad D^{\mathbb{A}}\psi = 0 .$$

(2.1)

The notation is as follows: First, $*$ is the metric's Hodge dual and $F_{\mathbb{A}}$ denotes the curvature 2-form of $\mathbb{A}$. The i-valued 1-form $\psi^\dagger \tau \psi$ is defined to be dual to the Clifford homomorphism $cl: T^*M \to End(\mathbb{S})$ that obeys $cl(a)cl(a') = -\langle a, a' \rangle \mathbb{I} - cl(*(a \wedge a'))$. This is to say that $\langle a, \psi^\dagger \tau \psi \rangle = \psi^\dagger cl(a) \psi$. Finally, $D^{\mathbb{A}}$ denotes the Dirac operator that is defined using the Levi-Civita connection on $T^*M$ and the connection $\mathbb{A}$ on $S_F$.

The group $C^\infty(M; S^1)$ acts on the space of solutions to (2.1) as follows: View $S^1$ as the unit circle in $\mathbb{C}$. Any given element $u \in C^\infty(M; S^1)$ sends a solution $\mathfrak{c} = (\mathbb{A}, \psi)$ to the solution $u \cdot \mathfrak{c} = (\mathbb{A} - 2u^{-1}du, u\psi)$. A solution to (2.1) is said to be *irreducible* when $\psi$ is not identically zero. The stabilizer in $C^\infty(M; S^1)$ of an irreducible solution is the constant function 1. Reducible solutions are those with $\psi = 0$ everywhere. Such solutions exist if and only if $det(\mathbb{S})$ has torsion first Chern class.

Let Conn($S_F$) denote the Frechet space of smooth, unitary connections on $S_F$. The equations in (2.1) can be viewed as asserting that ($\mathbb{A}$, $\psi$) is a critical point of a certain function on Conn($S_F$) $\times$ $C^\infty$(M; $\mathbb{S}$). This function is denoted by $\mathfrak{a}$ and it is given by:

$$\mathfrak{a} = \tfrac{1}{2} \mathfrak{cs} - r\mathfrak{e} - \mathfrak{e}_\mu + r \int_M \psi^\dagger D_\mathbb{A} \psi ,$$

(2.2)

where $\mathfrak{cs}$, $\mathfrak{e}$ and $\mathfrak{e}_\mu$ are defined as follows: Fix a flat connection, $\mathbb{A}_\mathbb{S}$, on $S_F$ and write $\mathbb{A} = \mathbb{A}_\mathbb{S} + 2\hat{a}$. Then

$$\mathfrak{cs}(\mathbb{A}) = -\int_M \hat{a} \wedge d\hat{a}$$

(2.3)

Meanwhile

$$\mathfrak{e} = i \int_M \upsilon \wedge d\hat{a}$$

(2.4)

where $\upsilon$ is such that $d\upsilon = \omega$. Finally, $\mathfrak{e}_\mu$ is defined using (2.4) with $\mu$ replacing $\upsilon$. A third function also plays a distinguished role here, this

$$\mathrm{E} = i \int_M \omega \wedge *d\hat{a} .$$

(2.5)

The next proposition relates the Seiberg-Witten equations to the dynamics. To set the notation, remark that Clifford multiplication by the 1-form $*\omega$ defines a splitting of $\mathbb{S}$ as

$$\mathbb{S} = E \oplus EK^{-1} ,$$

(2.6)

where $K^{-1}$ is isomorphic as an oriented real bundle to the kernel in TM of the 1-form $*\omega$ with the orientation given by $\omega$. A section $\psi$ of $\mathbb{S}$ is writen as $\psi = (\alpha, \beta)$ with respect to this splitting.

**Proposition 2.1**: *Suppose that $\{r_n\}_{n=1,2...}$ is an unbounded sequence of real numbers such that for each n there is a solution, ($\mathbb{A}_n$, $\psi_n$), to the $r = r_n$ version of (2.1). Write $\psi_n$ with respect to the splitting in (2.6) as $\psi_n = (\alpha_n, \beta_n)$. Suppose that $\{\sup_M (1 - |\psi_n|)\}$ is bounded away from zero.*
- *If the sequence $\{\mathrm{E}(\mathbb{A}_n)\}_{n=1,2,...}$ has a convergent subsequence, then there is a subsequence of $\{(r_n, (\mathbb{A}_n, (\alpha_n, \beta_n))\}_{n=1,2,...}$ such that the corresponding sequence of* *measures $\sigma_n = r_n(1 - |\alpha_n|^2)^2$ dvol converges weakly to a measure with support on a non-empty union of closed, integral curves of* v.

- *If the sequence $\{E(\mathbb{A}_n)\}_{n=1,2,...}$ has no covergent subsequence, then the corresponding sequence of measures $\{\sigma_n = E(\mathbb{A}_n)^{-1} r_n(1 - |\alpha_n|^2)$ dvol$\}_{n=1,2,...}$ converges weakly to a non-trivial measure that is invariant with respect to the 1-parameter group of diffeomorphisms that is generated by* v.

This proposition is proved in the next section.

### 3. Proof of Proposition 2.1

Some preliminary results are required. The proofs are much the same as those that appear in [T1]. In what follows, $\nabla$ denotes the covariant derivatives on the sections of the summands E and $EK^{-1}$ of $\mathbb{S}$ that are induced by the given connection $\mathbb{A}$ on $S_F$.

**Lemma 3.1**: *There exist constants $\kappa_1, \kappa_2 > 1$ such that if $r \geq 1$ and if $(\mathbb{A}, \psi = (\alpha, \beta))$ is a solution to (2.1), then*

- $|\alpha| \leq 1 + \kappa_1 r^{-1}$.
- $|\beta|^2 \leq \kappa_1 r^{-1}(1 - |\alpha|^2) + \kappa_2 r^{-2}$.

*In addition, given $k \geq 1$, there exists an r and $(A, \psi)$ independent constant $\kappa_k$ such that*

- $|\nabla^k \alpha|^2 + r|\nabla^k \beta|^2 \leq \kappa_k r^k$

*Finally, there exists $u \in C^{\infty}(M; U(1))$ such that $\mathbb{A} - 2u^{-1}du = \mathbb{A}_\mathbb{S} + 2\hat{a}$ where*

- $|\hat{a}| \leq \kappa r^{2/3}(|E|^{1/3} + 1)$

**Proof of Lemma 3.1**: The equation $D_\mathbb{A} \psi = 0$ implies that $D_\mathbb{A}^2 \psi = 0$. The Weitzenboch formula for the latter reads

$$\nabla^{\dagger}\nabla\psi + \tfrac{1}{4} R \psi + \tfrac{1}{2} \mathrm{cl}(*F_\mathbb{A})\psi = 0 ,$$

(3.1)

where R denotes the metric's scalar curvature. Taking the inner product with $\psi$ gives the inequality

$$\tfrac{1}{2} d^{\dagger}d|\psi|^2 + |\nabla\psi|^2 + r|\psi|^2(|\psi|^2 - 1) \leq c_0 |\psi|^2 ,$$

(3.2)

where $c_0 \geq 1$ is an r and $(\mathbb{A}, \psi)$ independent constant. This use of $c_0$ is the first illustration of the following convention used in the remainder of this paper: In all cases, $c_0$ denotes a constant that is at least 1 and independent of r and of any given pair of connection and section of $\mathbb{S}$. Its value can be assumed to increase from appearance to appearance.

An application of the maximum principle to this last equation finds $|\psi|^2 \leq 1 + c_0 r^{-1}$. To obtain estimates on $\alpha$ and $\beta$, introduce the projections $P_\pm: \mathbb{S} \to \mathbb{S}$ onto the respective E and $EK^{-1}$ summands. Apply these projections to the right hand side of (3.1) to see that

- $\nabla^\dagger \nabla \alpha + r(|\alpha|^2 + |\beta|^2 - 1)\alpha + \mathcal{R}_1 \cdot \nabla\alpha + \mathcal{R}_2 \cdot \nabla\beta + \mathcal{R}_3 \alpha + \mathcal{R}_4 \beta = 0$
- $\nabla^\dagger \nabla \beta + r(|\alpha|^2 + |\beta|^2 + 1)\beta + \mathcal{S}_1 \cdot \nabla\alpha + \mathcal{S}_2 \cdot \nabla\beta + \mathcal{S}_3 \alpha + \mathcal{S}_4 \beta = 0.$

(3.3)

Let $w = (1 - |\alpha|^2)$. These two equations imply

- $\tfrac{1}{2} d^\dagger dw + r|\alpha|^2 w - \tfrac{1}{2} |\nabla\alpha|^2 - r|\alpha|^2|\beta|^2 \geq -c_0(|\alpha|^2 + |\nabla\beta|^2 + |\beta|^2)$
- $\tfrac{1}{2} d^\dagger d|\beta|^2 + r|\alpha|^2|\beta|^2 + \tfrac{1}{2} r|\beta|^2 + |\nabla\beta|^2 \leq c_0 r^{-1}(|\nabla\alpha|^2 + |\alpha|^2).$

(3.4)

It follows as a consequence that there exists $c_0, c_1 \geq 1$, both independent of r and $(\mathbb{A}, \psi)$ such that $x = |\beta|^2 - c_0 r^{-1} w - c_1 r^{-2}$ obeys $d^\dagger dx + r|\alpha|^2 x \leq 0$. Thus, $x$ has no positive maximum. This proves the second assertion. To obtain the assertions about the derivatives of $\alpha$ and $\beta$, take Gaussian normal coordinates near any given point, then rescaling in these coordinates so the ball of radius $r^{-1/2}$ centered on the point has radius 1. Differentiate the rescaled versions of (3.3) and use the Green's function for the Laplacian in $\mathbb{R}^3$ to find in order $\kappa_1, \kappa_2, \ldots$, etc. The proof of the final assertion mimics that given for Lemma 2.4 in [T1].

The projections $P_\pm$ can be used to relate certain derivatives of $\alpha$ to derivatives of $\beta$. In particular, one can write the Dirac equation in 2-component form schematically as

$$\tfrac{1}{2} \hat{\nabla}_v \alpha + \bar{\partial}^\dagger \beta = 0 \quad and \quad \tfrac{1}{2} \hat{\nabla}_v \beta - \bar{\partial} \alpha = 0$$

(3.5)

where $\hat{\nabla}_v$ is the covariant derivative in the direction of v as defined using a unitary connection that is obtained from $\mathbb{A}$, $\omega$ and the Levi-Civita connection. Meanwhile, $\bar{\partial}$ denotes a version of the d-bar operator for derivatives orthogonal to v that is defined by this same data. In particular, near any given point, $\bar{\partial}$ can be written as $\tfrac{1}{2}(\hat{\nabla}_{e_1} + i \hat{\nabla}_{e_2})$ where $\{v, e_1, e_2\}$ define an orthonormal basis for TM.

*Proof of Proposition 2.1*: Consider a sequence $\{r_n, (\mathbb{A}_n, \psi_n = (\alpha_n, \beta_n))\}_{n=1,2,\ldots}$ as in the statement of Proposition 2.1. Suppose first that $\{E_n = E(\mathbb{A}_n)\}_{n=1,2,\ldots}$ has a convergent subsequence. Given Lemma 3.1, the arguments in Section 6 of [T1] that prove the latter's Theorem 2.1 can be employed with only cosmetic changes to prove what is asserted by the first bullet. The details of this are left to the reader.

To prove the second bullet, suppose that $\{E_n\}_{n=1,2,\ldots}$ has no convergent subsequence. It is a consequence of the first item in Lemma 3.1 that this sequence must be unbounded from above. This understood, pass to an increasing subsequence with each member greater than 1. Let $\{\sigma_n\}_{n=1,2,\ldots}$ denote the sequence of measures defined on $C^\infty(M)$ given by

$$f \to \sigma_n(f) = E_n^{-1} r_n \int_M f(1 - |\alpha_n|^2).$$

(3.6)

It follows from Lemma 3.1 using the top equation in (2.1) that $|\sigma_n(f)| \leq (1 + c_0 r^{-1/2})\|f\|_\infty$. As a consequence, the sequence $\{\sigma_n\}_{n=1,2,\ldots}$ is bounded in the dual to $C^0(M)$. Thus, it has weakly convergent subsequence with non-trivial limit; this because $|\sigma_n(1) - 1| \leq c_0 r^{-1/2}$. Let $\sigma_\infty$ denote the limit measure.

To prove that $\sigma_\infty$ is invariant, suppose for the moment that $(\mathbb{A}, \psi = (\alpha, \beta))$ is a given solution to some $r \geq 1$ version of (2.1) with $E(\mathbb{A}) \geq 1$. Let $f$ denote a $C^1$ function on $M$ with supremum norm 1. An integration by parts finds

$$r \int_M v(f)(1 - |\alpha|^2) = r \int_M f(\alpha \hat{\nabla}_v \bar\alpha + \bar\alpha \hat{\nabla}_v \alpha).$$

(3.7)

To continue, use the right most equation in (3.5) to write

$$r \int_M v(f)(1 - |\alpha|^2) = -r \int_M f(\alpha \overline{\bar\partial^\dagger \beta} + \bar\alpha \bar\partial^\dagger \beta).$$

(3.8)

Now integrate by parts to see that

$$\left| r \int_M v(f)(1 - |\alpha|^2) \right| \leq c_0 r \int_M |\hat\nabla \beta||\beta| + c_0 r \int_M |df||\beta|.$$

(3.9)

According to Lemma 3.1, the right hand side of (3.9) is bounded by

$$c_0 r \int_M (|f| + |df|)|\beta|.$$

(3.10)

To bound the latter, fix $\delta > 0$ and use the second item in Lemma 3.1 to bound (3.10) by

$$c_0 (1 + \|df\|_\infty)(\delta E + \delta^{-1}).$$

(3.11)

This understood, take $\delta = E^{-1/2}$ in (3.10) to see that

$$|r \int_M v(f) (1 - |\alpha|^2)| \leq c_0 (r^{-1/2} + E^{-1/2}) E (1 + \|df\|_\infty) .$$

(3.12)

With the preceding understood, take $(A, \psi) = (A_n, \psi_n)$ to see that

$$|\sigma_n(v(f))| \leq c_0 (r_n^{-1/2} + E_n^{-1/2})(1 + \|df\|_\infty) .$$

(3.13)

Thus, $\sigma_\infty(v(f)) = 0$ if $f$ is any given $C^1$ function. This implies the second bullet of the Proposition 2.1 since $C^1(M)$ is dense in $C^0(M)$.

## 4. Input for Proposition 2.1

Let $\mathfrak{c} = (\mathbb{A}, \psi)$ denote a pair of connection on $S_F$ and section of $\mathbb{S}$. For each $r \geq 1$, associate to $\mathfrak{c}$ the operator on $C^\infty(M; iT^*M \oplus \mathbb{S} \oplus i\mathbb{R})$ that sends a given section $(b, \eta, \phi)$ to the section whose respective $iT^*M$, $\mathbb{S}$ and $i\mathbb{R}$ components are

- $*db - d\phi - 2^{-1/2} r^{1/2} (\psi^\dagger \tau \eta + \eta^\dagger \tau \psi)$,
- $D^\mathbb{A} \eta + 2^{1/2} r^{1/2} (cl(b)\psi + \phi \psi)$,
- $*d*b - 2^{-1/2} r^{1/2} (\eta^\dagger \psi - \psi^\dagger \eta)$.

(4.1)

This operator is symmetric, and it has an unbounded, self-adjoint extension to the space of square integrable sections of $iT^*M \oplus \mathbb{S} \oplus i\mathbb{R}$. This extension has compact resolvent and so its spectrum is discrete with finite multiplicities and no accumulation points. This self-adjoint extension of (4.1) is denoted in what follows by $\mathfrak{L}_{\mathfrak{c}, r}$. Fix a section $\psi_\mathbb{S}$ of $\mathbb{S}$ so that the $r = 1$ and $(\mathbb{A}_\mathbb{S}, \psi_\mathbb{S})$ version of (4.1) has trivial kernel. First order perturbation theory can be used to find such a section.

Fix $r \geq 1$ and $\mathfrak{c} \in \text{Conn}(S_F) \times C(M; \mathbb{S})$ such that $\mathfrak{L}_{\mathfrak{c}, r}$ has trivial kernel. Associate to $r$ and $\mathfrak{c}$ an integer $f(\mathfrak{c}, r) \in \mathbb{Z}$ as follows: Fix any path, $s \to \mathfrak{c}(s)$, from $[1, r]$ into $\text{Conn}(S_F) \times C^\infty(M; \mathbb{S})$ with $\mathfrak{c}(1) = (\mathbb{A}_F, \psi_F)$ and $\mathfrak{c}(r) = \mathfrak{c}$. The integer $f(\mathfrak{c}, r)$ is the spectral flow for the path of operators, $\{\mathfrak{L}_{\mathfrak{c}(s), s}\}_{s \in [1, r]}$. See, e.g. [T4]. Note that $f(\mathfrak{c}, r) = f(u \cdot \mathfrak{c}, r)$ for any map $u: M \to S^1$.

The next proposition will be used to obtain input for the application of Proposition 2.1. To set the stage, remark that there exist, in all cases, Spin$^{\mathbb{C}}$ structures on M with torsion first Chern class. This follows by virtue of the fact that oriented 3-manifolds are spin manifolds and so the first Chern class of $K^{-1}$ is divisible by 2.

**Proposition 4.1**: *Suppose that the Spin$^{\mathbb{C}}$ structure has torsion first Chern class. There is an unbounded set $\Lambda \in \mathbb{Z}$, $\kappa > 0$ and, given $\varepsilon > 0$, integers $m \geq 3$ and $N \geq 1$, a 1-form $\mu$ with $C^m$ norm less than $\varepsilon$; all with the following properties: Fix an integer $f \in \Lambda$ with $f \leq N$. Then there exists*

- *a set $\{\rho_k\}_{k=0,1,\ldots} \subset [1, \infty)$ with no accumulation points.*
- *a continuous map $\mathfrak{a}_\Diamond: [\rho_0, \infty)$;*
- *for each index k, a smooth map $\hat{\mathfrak{c}}_k: [\rho_k, \rho_{k+1}] \to \mathrm{Conn}(S_F) \times C^\infty(\mathbb{S})$ such that $\hat{\mathfrak{c}}_k(r)$ is an irreducible solution to the $(r, \mu)$ version of (2.1) for each $r \in [\rho_k, \rho_{k+1}]$. Moreover,*
  a) *The operator $\mathfrak{L}_{\hat{\mathfrak{c}}(r),r}$ has trivial kernel.*
  b) *The spectral flow $f(\hat{\mathfrak{c}}_k(r), r)$ is equal to $f$.*
  c) *Write $\hat{\mathfrak{c}}_k(r) = (\mathbb{A}, \psi)$. Then $\sup_M (1 - |\psi|) \geq \kappa > 0$.*

*Finally, the function $r \to \mathfrak{a}(\hat{\mathfrak{c}}_k(r))$ is the restriction to $[\rho_k, \rho_{k+1}]$ of the function $\mathfrak{a}_\Diamond$.*

*Proof of Proposition 4.1*: Except for c) of the third bullet, the arguments are identical to those used to prove Propositions 4.2 in [T1]. To elaborate, Kronheimer and Mrowka introduce in Chapter 3 of [KM] a version of Seiberg-Witten Floer homology that they denote by $\widehat{H}_*$. For a Spin$^{\mathbb{C}}$ structure with torsion first Chern class, this is a relative $\mathbb{Z}$-graded $\mathbb{Z}$-module which is finitely generated in each degree. As explained in Chapter 35.1 of [KM], it is non-zero for an infinite set of degrees, a set that is bounded from above but not from below. The 1-form $\mu$ for use in Proposition 4.1 is chosen as in Section 3 of [T1] so that the irreducible solutions to (2.1) for $r \notin \{\rho_k\}$ can be identified with the a subset of the generators of a chain complex with differential whose cohomology is $\widehat{H}_*$. The relative degree between any two such generators is minus the difference between the corresponding values of the spectral flow function. The remaining generators are associated to the reducible solutions of (2.1). For $r \gg 1$, these all have very negative degree: It is a consequence of Proposition 5.5 in [T1] that the degrees of the latter are bounded from above by $-c_0^{-1} r^2 s_v$. Granted all of this, fix a degree, $f$, where $\widehat{H}_*$ has a non-zero class, choose a sufficiently generic $\mu$ as done in Section 3 of [T1], and then for all $r \notin \{\rho_k\}$ large, this class is represented by an irreducible solution to (2.1) such that the corresponding version of (4.1) has trivial kernel. The construction of the function $\mathfrak{a}_\Diamond$ and the maps $\{\hat{\mathfrak{c}}_k\}_{k=0,1,\ldots}$ is done by mimicking what is done in Sections 3 of [T1] and the proof of Proposition 4.2 of [T1]. The argument for c) of the third bullet is the same as that used to prove Lemma 5.4 in [T1].

Invoke Proposition 4.1 for a given small $\varepsilon$ and large $m$ and $N$, and integer $n \in \Lambda$. Define the piecewise differentiable (but perhaps discontinuous) functions $\mathfrak{cs}$, $\mathfrak{e}$, and $\mathrm{E}$ on $[\rho_0, \infty)$ by setting each on any given $[\rho_k, \rho_{k+1}]$ to be the composition of its namesake as defined in (2.3), (2.4) and (2.5)) with $\hat{c}_k(\cdot)$.

**Lemma 4.2**: *There is an unbounded sequence $\{r_n\}_{n=1,2,\ldots} \subset [\rho_0, \infty)$ such that either*
- *The corresponding sequence $\{E(r_n)\}_{n=1,2,\ldots}$ is bounded as $r \to \infty$.*
- *The sequence $\{E(r_n)\}_{n=1,2,\ldots}$ is divergent and the sequence $\{\max(\mathfrak{e}(r_n)/E(r_n), 0)\}_{n=1,2,\ldots}$ converges to zero.*

*Proof of Lemma 4.2*: There are various cases to consider.

Case 1: *There exists an unbounded set of value for $r$ where $\mathfrak{e} \leq 1$. In addition:*
    a) *The corresponding sequence of values for $E$ has a bounded subsequence.*
    b) *The corresponding sequence of values for $E$ has no bounded subsequence.*

Case 2: *There exists $r_*$ such that $\mathfrak{e} \geq 1$ when $r \geq r_*$. In addition, there exists $\delta \in (0, \frac{1}{8})$ and $z \in [0, \frac{1}{256}]$ such that either:*
    a) *There exists an unbounded set of $r$ in $[r_*, \infty)$ such that $E \geq \delta r^z \mathfrak{e}$.*
    b) $E \leq \delta r^z \mathfrak{e}$ *for all $r \geq r_*$.*

(4.2)

The conclusions of the lemma follow immediately if Cases 1a,b or 2a are relevant. The only problematic case is 2b. The discussion of the latter case requires the following lemma:

**Lemma 4.3**: *There exists $\kappa \geq 1$ with the following significance: Suppose that $r \geq 1$ and that $(\mathbb{A}, \psi)$ is a solution to (2.1) for a $Spin^{\mathbb{C}}$ structure with torsion first Chern class. Then $\mathfrak{cs} \leq \kappa\, r^{2/3}(|E|^{4/3} + \kappa)$.*

    To consider the implications of Case 2b, let $r \to \mathfrak{v}(r)$ denote $-2\mathfrak{a}_\lozenge(r)/r$. This function is continuous and piecewise differentiable. In particular,

$$\tfrac{d}{dr}\mathfrak{v} = -\mathfrak{v}/r + \mathfrak{e} = \mathfrak{cs}/r^2$$

(4.3)

on $(\rho_k, \rho_{k+1})$. Note as well that

$$\tfrac{d}{dr}(r\mathfrak{v}) = \mathfrak{e}\;;$$

(4.4)

and so the function $r \to r\,\mathfrak{v}(r)$ is increasing where $\mathfrak{e} > 0$.

To consider the implications of this for Case 2b, note that there exists $r_*$ such that $E \geq 1$ for all $r \geq r_*$. Granted this, there exists $\varepsilon \geq \frac{1}{256}$ such that

$$\mathfrak{cs} \leq c_0 r^{1-\varepsilon} E \quad \textit{for all} \ \ r \geq r_*.$$

(4.5)

Indeed, if such were not the case, Lemma 4.3 would give an infinite subsequence of r values with $E \geq c_0^{-1} r^{1-3\varepsilon}$, and hence with $\mathfrak{cs} \geq c_0^{-1} r^{2-4\varepsilon}$. This would imply, via Lemma 5.3 Proposition 5.5 in [T1] that the spectral flow functions for the various $\hat{\mathfrak{c}}_k$ were not identical.

Given (4.5), it then follows in Case 2b that $\mathfrak{cs} \leq \delta r^{1-z} \mathfrak{e}$ and so $\mathfrak{v} \geq (2 - \delta) \mathfrak{e}$. This implies that $\mathfrak{cs} \leq \delta r^{1-z} \mathfrak{v}$. Use the latter bound in (4.3) to see that

$$\tfrac{d}{dr} \mathfrak{v} \leq \delta r^{1-z} \mathfrak{v}.$$

(4.6)

This last bound integrates to give $\ln(\mathfrak{v}(r)/\mathfrak{v}(r_*)) \leq \delta z^{-1} r_*^{-z}$. As a consequence, the function $r \to \mathfrak{v}(r)$ is bounded as $r \to \infty$ if $z > 0$. Note that this implies that $\mathfrak{e}$ is also bounded because $\mathfrak{v} \leq (2 + \delta)\mathfrak{e}$. Thus, Case 2b requires that $\mathfrak{e}$ is bounded if $z > 0$.

A bound for $\mathfrak{e}$ also follows if $z = 0$ provided that $\delta$ is small. To see that such is the case, use Lemma 4.3s bound $\mathfrak{cs} \leq c_0 r^{2/3} E^{4/3}$ to see that $\mathfrak{cs} \leq c_0 r^{2/3} \delta^{4/3} \mathfrak{v}^{4/3}$. This with (4.3) gives

$$-3 \tfrac{d}{dr} \mathfrak{v}^{-1/3} \leq c_0 \delta^{4/3} r^{-4/3},$$

(4.7)

which integrates to yield

$$3\mathfrak{v}(r_*)^{-1/3} - 3\mathfrak{v}(r)^{-1/3} \leq c_0 \delta^{4/3} r_*^{-1/3}.$$

(4.8)

This gives a bound on $\mathfrak{v}(r)$. Indeed, because $E < c_0 r$, and $\mathfrak{e} \leq c_0 r$, it follows that $\mathfrak{v} \leq c_0 r$, and so $\mathfrak{v}^{-1/3} \geq c_0^{-1} r^{-1/3}$. Thus, if $\delta$ is less than a positive constant given by $r_*$, then $\mathfrak{v}(r)$ is bounded; and thus so is $\mathfrak{e}$.

The conclusions of Lemma 4.2 are seen to hold in Case 2b because $E$ is not bounded and $\mathfrak{e}$.

*Proof of Lemma 4.3*: It follows from Lemma 3.1 and the left most equation in (2.1) that

$$\int_M |F_{\mathbb{A}}| \leq c_0 (E + 1)$$

(4.9)

when $(\mathbb{A}, \psi)$ is a solution to (2.1). In the case when $\det(\mathbb{S})$ is torsion, the function $\mathfrak{cs}$ is gauge invariant. This understood, fix a gauge such that $\mathbb{A} = \mathbb{A}_S + 2\hat{a}$ where $\mathbb{A}_S$ is flat and where â obeys $d^\dagger \hat{a} = 0$ and is such that the $L^2$ norm of the projection of â into the space of harmonic one forms is bounded by $c_0$. According to the final item in Lemma 3.1, the norm of the 1-from â s bounded pointwise by $c_0 (r^{2/3} |E|^{1/3} + 1)$. This and (4.9) imply the assertion of Lemma 4.3.

## 5. Proof of Theorems 1 and 2

To prove Theorem 1, invoke Proposition 4.1 for some small $\varepsilon$, large m and N, and some $f \in \Lambda$ with $f \leq N$. Let $\{r_n\}_{n=1,2,...}$ denote the sequence that is supplied by Lemma 4.2. Each $r_n$ is contained in some interval $[\rho_{k(n)}, \rho_{k(n)+1}]$. Use $\mathfrak{c}_n$ to denote $\hat{\mathfrak{c}}_{k(n)}(r_n)$. Now apply Proposition 2.1 to the sequence $\{\mathfrak{c}_n\}_{n=1,2,...}$. To see what results, it proves useful to write $\upsilon = q*\omega + \mathfrak{p}$ where $\mathfrak{p}$ is annihilated by v and where q is a smooth function. Note that the integral of q over M is $s_\upsilon$ and so non-zero. Suppose that $(\mathbb{A}, \psi = (\alpha, \beta))$ is a solution to the $r \gg 1$ version of (2.1). It follows from the left most identity in (2.1) and Lemma 3.1 that

$$\mathfrak{e} = r \int_M q(1 - |\alpha|^2) + \mathfrak{r},$$

(5.1)

where $|\mathfrak{r}| \leq c_0 E^{1/2}$. This understood, it follows that the subsequence given in Proposition 2.1 can be chosen so that the limit measure, $\sigma_\infty$, obeys $\sigma_\infty(q) \leq 0$. This implies that $\sigma_\infty$ is not integration with respect to the volume form on M because the latter pairs with q to give a positive number, $s_\upsilon$. Thus the dynamics is not uniquely ergodic.

To prove Theorem 2, choose write the 1-form a as $q*\omega + \mathfrak{p}$ where q is again a function M and where $\mathfrak{p}$ is is annihilated by v. In this case, $q \geq 0$. Fix $\delta > 0$ and let $M_\delta \subset M$ denote the set of points where $q \geq \delta$. It now follows from Lemma 5.1 that

$$\mathfrak{e} \geq r\delta \int_{M_\delta} (1 - |\alpha|^2).$$

(5.2)

This understood, invoke Proposition 4.1 for some small $\varepsilon$, large m and N, and some $f \in \Lambda$ with $f \leq N$. Let $\{r_n\}_{n=1,2,...}$ denote the sequence that is supplied by Lemma 4.2. If the sequence $\{E(r_n)\}_{n=1,2,...}$ is not bounded, then $\lim_{n \to \infty} \mathfrak{e}(r_n)/E(r_n) \to 0$. To see what this implies, introduce for each n, the solution $\mathfrak{c}_n$ as in the preceding paragraph. Write $\mathfrak{c}_n$ as $(\mathbb{A}_n, \psi_n = (\alpha_n, \beta_n))$. Then it must be the case that for any fixed $\delta > 0$,

$$\lim_{n \to \infty} E(\mathfrak{c}_n)^{-1} \int_{M_\delta} (1 - |\alpha_n|^2) = 0.$$

(5.3)

This implies that the measure given by Proposition 2.1 is supported on the set where the function q = a(v) is zero.

## 6. An example: The horocycle flow

The author learned of this example from Curt McMullen. Let Σ denote a compact surface of genus greater than 1 and with a hyperbolic metric. Let M → Σ denote its unit tangent bundle. The 3-manifold M has a volume preserving, 1-parameter family of diffeomorphisms that is called the 'horocycle flow'. This flow is generated by a vector field, $v_H$, with the following two properties: First, v is uniquely ergodic. Second, ω = Ω($v_H$, ·) is exact. According to Theorem 1, the corresponding version of (1.1) must vanish. The purpose of this section is to verify directly that such is the case. This follows from

**Lemma 6.1**: *The form ω for the vector field $v_H$ has an anti-derivative, υ, which is nowhere zero and is such that υ ∧ ω = 0. In particular, the kernel of υ defines a foliation of M. Moreover, both the vector field $v_H$ and the generator of the geodesic flow on M are annihilated by υ; thus they span the tangent space to the leaves of the foliation.*

***Proof of Lemma 6.1***: Let $\mathbb{H} \subset \mathbb{C}$ to denote the upper half plane, thus the set {z = x + iy ∈ $\mathbb{C}$: y > 0}. The hyperbolic metric on $\mathbb{H}$ has line element

$$ds^2 = y^{-2}(dx^2 + dy^2).$$

(6.1)

Use this metric to define the unit tangent bundle $S\mathbb{H}$. This bundle has coordinates (z, θ) with z ∈ $\mathbb{H}$ and with θ defined by writing a vector $u \frac{\partial}{\partial z}$ as $|u| e^{i\theta} \frac{\partial}{\partial z}$. Note that (z, u) is in the unit tangent bundle $S\mathbb{H}$ when |u| = y. The volume form on $S\mathbb{H}$ is Ω = $y^{-2}$ dxdydθ.

Define a map φ: Sl(2; ℝ) → $T^{1,0}\mathbb{H}$ as follows. Write an element in $T^{1,0}\mathbb{H}$ as a pair (z, $u\frac{\partial}{\partial z}$) and write an element $\mathbb{A}$ ∈ Sl(2; ℝ) as

$$\mathbb{A} = \begin{pmatrix} a & b \\ c & d \end{pmatrix} \in SL(2; \mathbb{R}).$$

(6.2)

Here, ad - bc = 1. let ς = b + ia and σ = d + ic. Then φ(A) is the point with coordinates

$$z = z(\mathbb{A}) = \varsigma/\sigma \in \mathbb{H} \quad and \quad u = i\sigma^{-2}.$$

(6.3)

Such a point lies in $S\mathbb{H}$ when $y|\sigma|^2 = 1$. This map gives the well known identification between $S\mathbb{H}$ and $SL(2; \mathbb{R})$.

The horocycle flow on $T^{1,0}\mathbb{H}$ is induced by the flow $\mathbb{R} \times Sl(2; \mathbb{R}) \to Sl(2; \mathbb{R})$ that sends

$$(t, \mathbb{A}) \to \begin{pmatrix} a & b+ta \\ c & d+tc \end{pmatrix}.$$

(6.4)

The generator of this flow is the left-invariant vector field $v_H$ given by the element

$$\sigma_+ = \begin{pmatrix} 0 & 1 \\ 0 & 0 \end{pmatrix}$$

(6.5)

in the Lie algebra of $SL(2; \mathbb{R})$. Together with $e_+$, the matrices $e_- = e_+^T$ and

$$\sigma_3 = [\sigma_+, \sigma_-] = \begin{pmatrix} 1 & 0 \\ 0 & -1 \end{pmatrix}$$

(6.6)

span the Lie algebra of $SL(2; \mathbb{R})$. The left-invariant vector field, $v_G$, that corresponds to $\sigma_3$ generates the geodesic flow on $S\mathbb{H}$; this the image in $S\mathbb{H}$ of the map from $\mathbb{R} \times Sl(2; \mathbb{R})$ to $Sl(2; \mathbb{R})$ that sends $(t, \mathbb{Z})$ to

$$\begin{pmatrix} ae^t & be^{-t} \\ ce^t & de^{-t} \end{pmatrix}.$$

(6.7)

The commutation relation $[\sigma_3, \sigma_+] = 2\sigma_+$ imply that $v_H$ and $v_G$ together span an integrable subbundle, $F_H \subset T(S\mathbb{H})$. The latter is therefore tangent to the leaves of a foliation.

Let $\upsilon^+$, $\upsilon^3$ and $\upsilon^-$ denote the corresponding, dual left-invariant vector fields to the vector fields $v_H$, $v_G$ and the left invariant vector field generated by $\sigma_-$. The volume form, $\Omega$, appears on $SL(2; \mathbb{R})$ as the left-invariant form $\Omega = 2\, \upsilon^+ \wedge \upsilon^3 \wedge \upsilon^-$. Thus $\Omega(v_H, \cdot) = \omega = 2\upsilon^3 \wedge \upsilon^-$. The as yet unstated commutation relation $[\sigma_3, \sigma_-] = -2\sigma_-$ implies that $\omega$ can be written as $d\upsilon^-$; this verifies that $\omega$ is exact. Note that the subbundle $F_H$ is the kernel of $\upsilon^-$. In any event, $\omega = 2\,\upsilon^3 \wedge \upsilon^-$ has the anti-derivative $\upsilon = \upsilon^-$ with the property that $\upsilon \wedge \omega = 0$.

As all of the forms and vector fields discussed above are left-invariant, the preceding discussion applies verbatim to quotients of $SL(2; \mathbb{R})$ by the action from the left by discrete subgroups. As a consequence, what is said above applies to the unit tangent space of any compact hyperbolic surface. The flow generated by $v_H$ on such a unit tangent space is uniquely ergodic, it is exact, and $s_{v_H} = 0$ since $\Omega(v_H, \cdot) = \omega$ has an anti-derivative, $\upsilon$, with $\upsilon \wedge \omega = 0$.